\documentclass[twocolumn]{article}
\usepackage{graphicx}
\usepackage{hyperref}
\usepackage{xcolor}
\usepackage{comment}
\usepackage{amsmath, amssymb}
\usepackage{tikz}
\usepackage{siunitx}

\title{\texttt{fkcompute}, an efficient $F_K$ invariant calculator}
\date{July 2026}

\author{
Paul Orland,
Davide Passaro,
Lara San Mart\'in Su\'arez,
Toby Saunders-A'Court, \\
and Josef Svoboda \\[1ex]
\small{\textit{California Institute of Technology, Pasadena, CA 91125, USA}}
\thanks{Contact: \texttt{dpassaro@caltech.edu}, \texttt{svo@caltech.edu}.}
}

\makeatletter
\renewcommand\@fnsymbol[1]{}
\makeatother

\begin{document}

\maketitle

\begin{abstract}
We introduce \texttt{fkcompute}, an open-source package for computing the Gukov--Manolescu invariant of links from a braid presentation.
\texttt{fkcompute} implements Park's inverted state sum through a three-phase pipeline.
First, a search is performed for a suitable braid presentation and for an additional inversion data on the braid.
Then, the state space of the inverted sum is encoded as a polytope, bounded by the associated linear constraint system.
Finally, the invariant is constructed by multiplication of R-matrices associated to the states.
Benchmarks show that prime knots up to 12 crossings, and prime links up to 10 crossings and of at most 3 components, are comfortably within reach.
As a result, \texttt{fkcompute} is used to compile the first public database of the Gukov--Manolescu invariant.
The package is available as a Python library, a command-line tool, and a Mathematica paclet.
\end{abstract}

\tableofcontents

\section{Introduction}\label{sec:introduction}

Modern topology increasingly relies on the computation of invariants to distinguish, classify, and analyze geometric and algebraic structures.
For knots, links, and three-manifolds, these invariants are used not only to distinguish examples, but also to test conjectures, build tables, and guide new theoretical work.

This computational role has grown as the invariants themselves have become more refined.
Leading invariants have grown in sophistication, with their computation correspondingly becoming more complex and resource-intensive.
Classical invariants such as the Alexander polynomial and the Jones polynomial are relatively inexpensive to compute.
By contrast, modern invariants such as colored Jones polynomials, Khovanov homology, knot Floer homology, and the Gukov--Manolescu series, involve larger algebraic structures, recursive formulas, or state sums whose cost grows quickly with the diagrammatic complexity of the knot involved.

At the same time, large datasets of invariants have become increasingly important.
Traditional mathematical investigation is now often complemented by data-driven methods, including techniques from data science and machine learning, which are beginning to play a substantive role in contemporary research.
This shift makes reproducibility and accessibility more urgent.
As the number of invariants grows, and as their computation depends on increasingly specialized algorithms, it becomes harder to distribute implementations, verify results, and collaborate on further development.
The problem is especially pronounced when the relevant software consists of complex, domain-specific code.

This trend has raised the demand for higher accuracy in computation and more comprehensive datasets.
If invariants are to be used for experimentation, classification, and comparison, then their values must be computed reliably and made available in a form that other researchers can inspect, reproduce, and extend.
Datasets are valuable, but a dataset without the code that produced it is difficult to audit.
For newer invariants, this point is especially important: different normalizations, expansion conventions, and truncation choices can lead to superficially different outputs.

In this note, we present \texttt{fkcompute}, a Python package which aims to fill the current gap for \emph{open} computational tools for the $F_K$ invariant.
The $F_K$ invariant emerged from the study of the $\widehat{Z}$ invariant, a $q$-series quantum invariant of closed 3-manifolds \cite{GPV,GPPV}.
From this perspective, the $F_K$ invariant can be interpreted as a $\widehat{Z}$ invariant for knot and link complements\footnote{
When referring to links, the $F_K$ invariant is sometimes denoted in literature as $F_L$. However, in this note, we will use the two notations interchangeably.
}.
More broadly, the $\widehat{Z}$ invariant itself was introduced as part of a program aimed at developing categorified invariants of 3-manifolds, analogous in spirit to how Khovanov homology categorifies the Jones polynomial of knots \cite{Kh}, and $F_K$ is critical in this framework.

The \texttt{fkcompute} package provides a suite of utilities and algorithms that implement the necessary algebraic and combinatorial procedures to construct the invariant directly from braid data.
The implementation supports homogeneous braids and, more generally, the class of \emph{nice} knots and links as defined by Park \cite{ParkThesis}, which conjecturally includes all fibered knots; thereby covering a broad and practically relevant set of examples.
We tested the software on \emph{nice} knots with up to 12 crossings and \emph{nice} links up to 10 crossings and 3 components, demonstrating both robustness and computational feasibility within this range.
Alongside \texttt{fkcompute} we publish the first publicly available database of Gukov--Manolescu series for knots and links \cite{topologyfyi}.

\section{The \texttt{fkcompute} package}\label{sec:fkcompute-package}
\texttt{fkcompute} is provided as a Python package, as a standalone command-line interface, and as a Mathematica paclet enabling users to compute the $F_K$ invariant directly from a braid presentation in a flexible and accessible manner.
The package is designed with extensibility as a central goal: its source code is publicly available and structured to facilitate modification, extension, and integration into future computational workflows or research projects.

Users may interact with \texttt{fkcompute} programmatically by importing it into Python scripts, enabling automated computation, batch processing, and integration with existing mathematical software pipelines.
Alternatively, \texttt{fkcompute} can be used as a command-line tool, providing a straightforward interface for direct computation from the terminal without requiring additional programming.
For users working in Mathematica, we also provide a Wolfram Language Paclet (\texttt{FkCompute}) which calls the Python package, allowing computations and symbolic expressions to be used directly in notebooks.

\texttt{fkcompute} is freely available for download through its public GitHub repository \cite{fkcompute}, providing open access to both the source code and accompanying documentation.
Precomputed $F_K$ tables used for validation and benchmarking are distributed separately from the core source repository \cite{topologyfyi}.

Internally, \texttt{fkcompute} leverages high-performance open-source libraries to ensure computational efficiency and reliability: polynomial arithmetic is handled using the FLINT library, a C\texttt{++} library optimized for fast and precise computation with polynomials and power series, while the HiGHS mixed-integer optimization solver \cite{huangfu2018parallelizing} is used for constraint validation, enabling robust and efficient resolution of the linear and combinatorial conditions that arise in the construction of the invariant.

The package implements the inverted state-sum construction whenever an acceptable inversion datum is found.
In particular, it covers homogeneous braids and many non-homogeneous braid examples.
Each calculation takes as input a braid word, a truncation degree for the $x$-variables, and optionally a prescribed inversion datum.
After completing the computation it returns a truncated $F_K(x;q)$ series for knots, or a $F_L(x_1,\dots,x_\ell;q)$ series for links (where $\ell$ denotes the number of link components), as well as the inversion datum used and other metadata regarding the computation.

\section{Mathematical background}\label{sec:mathematical-background}

In this Section we recall some of the mathematical background of the $F_K$ invariant, as well as some of its properties which serve as checks on the computations performed by \texttt{fkcompute}.

\subsection{The \texorpdfstring{$F_K$}{FK} and \texorpdfstring{$F_L$}{FL} invariants}\label{sec:fk-fl-invariants}
The $F_K$ invariant originates from the same program which produced the $\widehat{Z}$ invariants for closed 3-manifolds \cite{GPV,GPPV}.
$\widehat{Z}$ invariants are $q$-series quantum invariants, motivated by considerations in supersymmetric quantum field theory in 6D.
They are expected to recover the WRT or Witten--Reshetikhin--Turaev \cite{Witten, RT} invariants through radial limits at roots of unity.
The $\widehat{Z}$ invariants are also expected to admit a categorification \cite{crane1994}, similar to that of the Khovanov homology for knots \cite{Kh}.
Gukov and Manolescu proposed the $F_L$ invariant as an analogue of this construction for knot complements \cite{GM}.
The $F_L$ invariant is a $\ell+1$ series, where $\ell$ is the number of components of the link it is describing.

The manifold and link invariants $\widehat{Z}$ and $F_L$ are intimately related to each other.
Since closed 3-manifolds can be obtained by Dehn surgery on links in $S^3$, one expects the $\widehat{Z}$ invariants of the resulting closed manifolds to be recoverable from the corresponding link-complement series \cite{GM}.
More precisely, if $S^3_{p/r}(K)$ denotes the result of $p/r$ surgery on a knot $K$, the expected surgery relation has the form
\begin{multline}\label{eq:surgery-background}
  \widehat{Z}_a\left(S^3_{p/r}(K)\right)
  =\\ \varepsilon q^\delta
  \mathcal{L}^{(a)}_{p/r}
  \left[
    \left(x^{\frac{1}{2r}}-x^{-\frac{1}{2r}}\right)F_K(x,q)
  \right].
\end{multline}
Here $\mathcal{L}^{(a)}_{p/r}$ denotes a ``Laplace'' transform determined by the surgery coefficient $p/r$ and by the label $a$, which corresponds to a $\operatorname{Spin}^c$ structure on the surgered manifold \cite[Thm 1.2]{GM}, for which a fully generic expression has not been found.
The extension of formula~\eqref{eq:surgery-background} is one of the principal motivations for the explicit computation of $F_K$: it relates the two-variable invariant of a knot complement to one-variable $q$-series invariants of closed 3-manifolds.

The values of the $F_K$ series for knots are also constrained by its relationship with the colored Jones polynomial.
Let $J_{K,n}(q)$ denote the reduced $n$-colored Jones polynomial, normalized so that the unknot has value $1$.
The Melvin--Morton--Rozansky \cite{MM,RozMMR,RozMMRConj} expansion studies $J_{K,n}(q)$ in the large-color regime
\begin{equation}
    q=e^\hbar, \qquad x=q^n=e^{n\hbar},
\end{equation}
where $x$ is held fixed and $\hbar\to 0$ in the real line.
In this limit, one obtains an asymptotic expansion of the form
\begin{equation}\label{eq:mmr-background}
  J_{K,n}(e^\hbar)
  \sim
  \sum_{j\geq 0}
  \frac{P_{K,j}(x)}{\Delta_K(x)^{2j+1}}\hbar^j,
\end{equation}
where $\Delta_K(x)$ is the Alexander polynomial of $K$.
The leading term in \eqref{eq:mmr-background} is the classical Melvin--Morton term and is given by the inverse of the Alexander polynomial.

Gukov and Manolescu conjectured that $F_K$ is obtained by reorganizing this perturbative expansion into an integral two-variable $q$-series \cite{GM}.
Park consequently proved in \cite{Park21} that this conjecture holds for a class of \emph{nice} knots constructed using his $R$-matrix approach. Explicitly, we expect the following relation to hold when expanding the rational function in the right-hand side at $x=0$:
\begin{equation}\label{eq:FK-mmr-background}
  \frac{F_K(x,q=e^{\hbar})}{x^{\frac12}-x^{-\frac12}}=\sum_{j\geq 0}\frac{P_K(j;x)}{\Delta_{K}(x)^{2j+1}}\frac{\hbar^j}{j!}~.
\end{equation}

\subsection{Large-color \texorpdfstring{$R$-matrices}{R-matrices}}\label{sec:large-color-r-matrices}

The computation of $F_K$ used in \texttt{fkcompute} follows Park's large-color $R$-matrix construction \cite{Park21}.
It is formally analogous to the Reshetikhin--Turaev construction of the colored Jones polynomial: a braid is decorated by representation-theoretic data, each crossing represents the action of the $R$-matrix, and the braid is closed by taking a reduced quantum trace.

We first fix the notation used throughout Sections~\ref{sec:large-color-r-matrices} and~\ref{sec:state-sums}.
Let $\beta$ be a braid on $N$ strands whose closure is the knot or link under consideration.
We orient the braid from bottom to top.
Let $V_\beta$ denote the set of crossings of $\beta$, and let $E_\beta$ denote the set of oriented segments obtained by splitting the braid at its crossings and identifying the top and bottom endpoints according to the braid closure.
For a crossing $c\in V_\beta$, write
\[
  e_{\mathrm{BL}}(c),\quad e_{\mathrm{BR}}(c),\quad
  e_{\mathrm{TL}}(c),\quad e_{\mathrm{TR}}(c)
\]
for the adjacent bottom-left, bottom-right, top-left, and top-right segments.
A state is a function
\[
  s:E_\beta\to \mathbb Z.
\]
At a crossing $c$, we abbreviate
\[
  \begin{split}
  i=s(e_{\mathrm{BL}}(c)),\qquad
  j=s(e_{\mathrm{BR}}(c)),\\
  i'=s(e_{\mathrm{TL}}(c)),\qquad
  j'=s(e_{\mathrm{TR}}(c)).
  \end{split}
\]
The local state variables are required to satisfy the conservation law
\begin{equation}\label{eq:state-conservation}
  i+j=i'+j',
\end{equation}
enforced in the formulas below by the Kronecker delta \(\delta_{i+j,i'+j'}\).

To each component of the link, we assign a \emph{topological variable}.
Let $x$ be the topological variable assigned to the strand running from the bottom-left segment to the top-right segment, and let $y$ be the topological variable assigned to the strand running from the bottom-right segment to the top-left segment.
Suppose first that $c$ is a positive crossing, in which case $x$ is the variable assigned to the overstrand and $y$ to the understrand. With these conventions, the extended large-color $R$-matrix is
\begin{multline}\label{eq:positive-rmatrix-1}
\check{R}(x,y)^{i',j'}_{i,j}=
\delta_{i+j,i'+j'}\,
q^{\frac{j+j'+1}{2}}
 x^{-\frac{i'+j+1}{4}}
 y^{-\frac{3j'-i+1}{4}}
  \\
\times q^{jj'}\begin{bmatrix} i \\ i-j' \end{bmatrix}_q
\prod_{1\le r\le i-j'}\left(1-q^{j+r}y^{-1}\right),
\end{multline}
if \(i\ge j'\ge 0\) or \(0>i\ge j'\),
\begin{multline}\label{eq:positive-rmatrix-2}
\check{R}(x,y)^{i',j'}_{i,j}=
\delta_{i+j,i'+j'}\,
q^{\frac{j+j'+1}{2}}
 x^{-\frac{i'+j+1}{4}}
 y^{-\frac{3j'-i+1}{4}}
  \\
  \times q^{jj'}\begin{bmatrix} i \\ j' \end{bmatrix}_q
\prod_{0\le r\le j'-i-1}\left(1-q^{j-r}y^{-1}\right)^{-1},
\end{multline}
if \(j'\ge 0>i\), and \(\check{R}(x,y)^{i',j'}_{i,j}=0\) otherwise.

For negative crossings, where now $x$ is assigned to the understrand and $y$ to the overstrand, we have:
\begin{multline}\label{eq:negative-rmatrix-2}
\check{R}^{-1}(x,y)^{i',j'}_{i,j}=
\delta_{i+j,i'+j'}\,
q^{-\frac{i+i'+1}{2}}
x^{\frac{3i'-j+1}{4}}
y^{\frac{j'+i+1}{4}}
\\
\times q^{-ii'}
\begin{bmatrix} j \\ i' \end{bmatrix}_{q^{-1}}
\prod_{0\le r\le i'-j-1}\left(1-q^{-i+r}x\right)^{-1},
\end{multline}
if \(i'\ge 0>j\),
\begin{multline}\label{eq:negative-rmatrix-1}
\check{R}^{-1}(x,y)^{i',j'}_{i,j}=
\delta_{i+j,i'+j'}\,
q^{-\frac{i+i'+1}{2}}
x^{\frac{3i'-j+1}{4}}
y^{\frac{j'+i+1}{4}}
\\
\times q^{-ii'}
\begin{bmatrix} j \\ j-i' \end{bmatrix}_{q^{-1}}
\prod_{1\le r\le j-i'}\left(1-q^{-i-r}x\right),
\end{multline}
if \(j\ge i'\ge 0\) or \(0>j\ge i'\),
and
\(
  \check{R}^{-1}(x,y)^{i',j'}_{i,j}=0
\)
otherwise.
Here
\begin{equation}\label{eq:qbinomial-background}
  (a;q)_k=\prod_{r=0}^{k-1}(1-aq^r),\qquad
  \begin{bmatrix} n \\ k \end{bmatrix}_q
  =\frac{(q^n;q^{-1})_k}{(q;q)_k}
\end{equation}
for \(k\geq 0\), with the convention that an empty product is equal to $1$.

For links, the topological variables \(x\) and \(y\) in the local formula are the variables assigned to the two link components passing through the crossing.
For knots, there is only one variable, so we can always set \(x=y\).

\subsection{State sums}\label{sec:state-sums}

The large-color $R$-matrices from Section~\ref{sec:large-color-r-matrices} determine the contribution for each crossing to the state sum.
To obtain a well-defined series, one must also specify which integer values are allowed on each segment of the braid.
This choice is encoded by an \emph{inversion datum}, a map
\[
  \iota:E_\beta\longrightarrow \{+,-\}.
\]
We use the convention that a segment marked \(+\) is assigned a non-negative integer and a segment marked \(-\) is assigned a negative integer:
\[
  \iota(e)=+ \Rightarrow s(e)\geq 0,
  \qquad
  \iota(e)=- \Rightarrow s(e)<0.
\]
The sign assignment must be locally compatible with the non-zero cases of the extended \(R\)-matrices.
Using notation
\[
\arraycolsep=1.4pt
  \begin{matrix}
  \iota(e_{\mathrm{TL}}) & \iota(e_{\mathrm{TR}})\\
  \iota(e_{\mathrm{BL}}) & \iota(e_{\mathrm{BR}})
  \end{matrix},
\]
the allowed local patterns are
\begin{equation}\label{eq:positive-inversion-patterns}
\begin{gathered}
\arraycolsep=1.4pt
\begin{matrix}-&+\\ +&-\end{matrix},\qquad
\begin{matrix}+&-\\ -&+\end{matrix},\qquad
\begin{matrix}-&-\\ -&-\end{matrix},\qquad
\begin{matrix}-&+\\ -&+\end{matrix},\qquad
\begin{matrix}+&+\\ +&+\end{matrix}
\end{gathered}
\end{equation}
for a positive crossing, and
\begin{equation}\label{eq:negative-inversion-patterns}
\begin{gathered}
\arraycolsep=1.4pt
\begin{matrix}-&+\\ +&-\end{matrix},\qquad
\begin{matrix}+&-\\ -&+\end{matrix},\qquad
\begin{matrix}-&-\\ -&-\end{matrix},\qquad
\begin{matrix}+&-\\ +&-\end{matrix},\qquad
\begin{matrix}+&+\\ +&+\end{matrix}
\end{gathered}
\end{equation}
for a negative crossing.
Thus, at every crossing, the number of incoming \(-\)-marked segments equals the number of outgoing \(-\)-marked segments.

Geometrically, we interpret the inversion datum as the multicycle on the braid defined by the segments mapped to $-$, where we allow ``jumps'' in the multicycle from an overstrand to an understrand (see \cite{Park21} for more details) in the local patterns
\begin{equation*}
\begin{gathered}
\arraycolsep=1.4pt
\begin{matrix}-&+\\ -&+\end{matrix},\qquad
\begin{matrix}+&-\\ +&-\end{matrix}
\end{gathered}.
\end{equation*}

The role of the inversion datum is partly formal and partly computational.
Formally, it specifies on which segments the geometric series expansions implicit in the extended $R$-matrix entries \eqref{eq:positive-rmatrix-2} and \eqref{eq:negative-rmatrix-2} are inverted, i.e.\ whether the state variable on a segment ranges over the non-negative or the negative integers.
The topological meaning of this choice is not yet fully understood.
What is essential for the present paper is that an inversion datum gives a concrete finite-order procedure for computing the corresponding state sum.

We will call an inversion datum \emph{acceptable} when each state sum \(x\)-coefficient is a finite Laurent polynomial in \(q\).
More precisely, for each fixed \(x\)-degree, only finitely many states should contribute.
This condition is what makes the formal sum computable: the coefficient of any fixed monomial in the \(x\)-variables can then be obtained by an ordinary finite calculation.

A basic example of inversion data is provided by homogeneous braids.
A braid is homogeneous if, for each generator index \(k\), all occurrences of \(\sigma_k\) in the braid word have the same sign.   Following Park \cite{Park21}, we label the right-hand segments of each positive (resp. negative) crossing by \(+\) (resp. by \(-\)). Homogeneity ensures that the labels are consistent along the braid.
The segments on the leftmost side of the braid may be labeled entirely by \(+\) or entirely by \(-\); either choice gives the same construction.

For non-homogeneous braids, there is no such simple rule.
Instead, we must search for an inversion datum satisfying the local sign conditions above and the coefficient-wise finiteness condition. The existence of an acceptable inversion datum is the practical condition that makes the inverted state sum computable.
Knots and links admitting such data are called \emph{nice}.

Let \(\Omega(\iota)\) denote the set of states \(s:E_\beta\to\mathbb Z\) satisfying the sign bounds determined by \(\iota\), the conservation law \eqref{eq:state-conservation}, and the local inequalities required for the corresponding \(R\)-matrix entry to be nonzero.
For \(s\in\Omega(\iota)\), let
\[
  R_c(s)=
  \begin{cases}
  \check R(x,y)^{i',j'}_{i,j}, & \text{ if }c \text{ is  positive},\\
  \check R^{-1}(x,y)^{i',j'}_{i,j}, & \text{ if } c \text{ is negative}
  \end{cases}
\]
where $x$ is the topological variable attached to the incoming overstrand (resp. understrand) and \(y\) to the incoming understrand (resp. overstrand) of a positive (resp. negative) crossing \(c\).

We use the reduced trace convention.
Let \(b_1,\ldots,b_N\) be the bottom segments of the braid, identified with the corresponding top segments after closure.
The state associated to the bottom of the first strand (equivalently, leftmost strand, when the braid is oriented bottom to top) is set to
\begin{equation}\label{eq:ground-state}
    s_0(b_1)=
  \begin{cases}
  0, & \iota(b_1)=+,\\
  -1,& \iota(b_1)=-.
  \end{cases}
\end{equation}
Define the subset of states
\[
  \Omega_0(\iota)=\{s\in\Omega(\iota):s(b_1)=s_0(b_1)\}.
\]
For each \(b_r\), \(r=2,\ldots,N\), define
\[
  \mu_r(s)=x_{\alpha(r)}^{-\frac12}q^{-\frac12-s(b_r)},
\]
where \(\alpha(r)\) is the link component containing \(b_r\).
In the knot case, this is simply \(x^{-\frac12}q^{-\frac12-s(b_r)}\).

For every state \(s\in\Omega_0(\iota)\), let
\begin{equation}\label{eq:state-contribution}
  P_\iota(s)
  =
  \left(\prod_{r=2}^{N} \mu_r(s)\right)
  \left(\prod_{c\in V_\beta} R_c(s)\right).
\end{equation}
Then, the inverted state sum is defined to be
\begin{equation}\label{eq:inverted-state-sum}
  Z^{\mathrm{inv}}_\iota(\beta)
  =
  (-1)^{\mathbf{s}(\iota)}
  \sum_{s\in\Omega_0(\iota)} P_\iota(s).
\end{equation}
Here, \(\mathbf{s}(\iota)\) is the number of closed components of the simple multicycle determined by the \(-\)-marked segments.
The open component containing $b_1$ is not counted.

The corresponding \(F_K\) series is
\begin{equation}\label{eq:FK-from-Zinv}
  F_K(x_1,\dots,x_\ell,q)=(x_1^{\frac12}-x_1^{-\frac12})Z^{\mathrm{inv}}_\iota(\beta).
\end{equation}
where $x_1=x_{\alpha(1)}$ is the variable associated to the component which contains $b_1$.

\section{Algorithm}\label{sec:algorithm}

The implementation in \texttt{fkcompute} evaluates the inverted state sum \eqref{eq:inverted-state-sum} by a three-phase pipeline.
Phase~1 enumerates possible multicycles, which encode candidate \emph{sign assignments}.
Phase~2 turns the definition of $\Omega_0(\iota)$ into a reduced integer linear constraint system, eliminates dependent state variables, performs the boundedness check needed for a coefficient-wise finite truncation to the requested $x$-degree, and exports the data in a compact ``ILP'' description.
Taken together, the simplified linear constraints define a polyhedron in the state variable space.
Phase~3 (a compiled C\texttt{++} binary) enumerates the integral points in this polyhedron and evaluates the $R$-matrices using FLINT-backed polynomial arithmetic.

\subsection{Topological Preprocessing and State Coordinates}\label{sec:topological-preprocessing}

Internally, a braid word is represented as a list \(\beta=[g_1,\dots,g_n]\) with each \(g_t\in\{\pm1,\dots,\pm N-1\}\).
The sign of \(g_t\) is the crossing sign and \(|g_t|=k\) means that, at height \(t\), strands \(k-1\) and \(k\) cross.
The Python class \texttt{BraidTopology}
realizes the closure of~\(\beta\) as a lattice diagram with coordinates \((i,t)\) where \(i\in\{0,\dots,N-1\}\) is the strand index and \(t\in\{0,\dots,n\}\) is the level between crossings.
For a crossing at height~\(t\) with \(|g_t|=k\), the four corners are
\begin{equation}\label{eq:corner-coordinates}
  (k-1,t),\ (k,t),\ (k-1,t+1),\ (k,t+1),
\end{equation}
which correspond to bottom-left, bottom-right, top-left, and top-right, respectively, for a braid oriented bottom to top.

Not every lattice position \((i,t)\) carries an independent state variable: if strand~\(i\) does not participate in the crossing at height~\(t\), then the integer label is transported unchanged through that layer.
Accordingly, \texttt{BraidTopology} maintains a map \(\mathrm{state}(i,t)\) which identifies \((i,t)\) with the canonical state variable representing the segment of the diagram containing that position.
Equivalently, the set of state variables is the set of segments obtained by cutting the braid diagram at the crossings, as in Section~\ref{sec:state-sums}.
The closure of the braid induces periodic identification of state variables
\begin{equation}\label{eq:periodicity}
  \mathrm{state}(i,n)\equiv \mathrm{state}(i,0),\qquad i=0,\dots,N-1,
\end{equation}
corresponding to identifications of top and bottom endpoints of the braid.
The component structure of the closure is computed by traversing the lattice diagram; the resulting component indices are used throughout to decide which of the topological variables \(\{x_j\}_{j=1}^\ell\) is attached to a given segment. Therefore, if $i<j$, then there is at least one strand of the component with topological variable $x_i$ that precedes the strands of the component with topological variable $x_j$.

\begin{figure*}[t]
    \centering
    \includegraphics[width=0.75\linewidth]{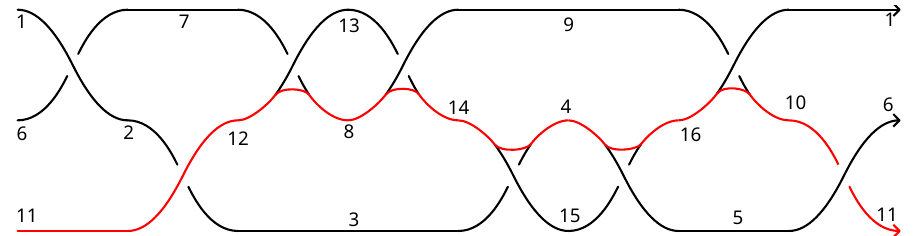}
    \caption{A sign assignment on the braid representative $[1,-2,-1,-1,2,2,-1,-2]$ of $8_{20}$, encoded as the multicycle associated to the permutation $\sigma = (11, 12, 8, 14, 4, 16, 10)$. The braid is oriented left to right.
    }
    \label{fig:8_20}
\end{figure*}

\subsection{Phase 1: Enumerating Multicycle Candidates}\label{sec:phase-1}

Phase~1 is implemented via
\texttt{find\_sign\_assignment}.
The braid is first classified as homogeneous or non-homogeneous. Each candidate determines a tentative sign assignment; local compatibility checks are performed at every crossing, and one of the four crossing types, denoted by \(R1,R2,R3,R4\), is assigned to each crossing,
corresponding exactly to the four nonzero cases in
\eqref{eq:positive-rmatrix-1}--\eqref{eq:negative-rmatrix-2} in that order.

\paragraph{Homogeneous braids}
If \(\beta\) is homogeneous, \texttt{fkcompute} constructs the inversion datum
deterministically using Park's rule, as described in Section~\ref{sec:state-sums}.

\paragraph{Non-homogeneous braids (search over multicycle candidates)}
For non-homogeneous braids, \texttt{fkcompute} enumerates representative \emph{multicycle candidates}.
The search is performed over braid variants obtained by cyclic rotation of the
braid word and, optionally, a horizontal mirror.
These transformations do not change the isotopy type of the closure but can dramatically change the outcome of the sign search.

The search begins by assigning a label to each segment of the braid, in the same order that they are encountered, walking the braid starting at its bottom left-most strand. If a link has more than one component, the label assignment is performed on every component of the link following the order of appearance in the braid strands, which coincides with the subscript of their topological variables. This gives labels in \(\{1,\dots,2n\}\) for a braid of $n$ crossings, and additional labels keep track of the component cycle.
As a result, we obtain four labels for every crossing, as assigned to its adjacent segments. The routine then records, for each crossing, the pair of labels associated to the incoming overstrand and understrand.

Any multicycle or inversion datum on a braid can be represented by a permutation on the labels \(\{1,\dots,2n\}\). Likewise, given a permutation $\sigma$, an associated multicycle is formed from the set of deranged points, by connecting the segment $i$ to $\sigma(i)$, while leaving the fixed points untouched (c.f. Figure \ref{fig:8_20}).

Using the equivalence above, all possible permutations associated to inversion data are generated. Given a crossing with incoming labels \((a,b)\), the value of the incoming overstrand label \(a\) may be left fixed, continued along the overstrand to \(a+1\), or it may ``jump'' to the outgoing understrand \(b+1\).
The value of the understrand label \(b\) may be left fixed or continued along the understrand to \(b+1\). Across all crossings, the chosen targets are required to be distinct, and the resulting output assembles into disjoint directed cycles in the closed braid diagram.
The induced sign assignment is read off by marking exactly the deranged points:
\begin{equation}\label{eq:continuation-to-signs}
  \iota(e)=\begin{cases}
    +1, & \text{if }\sigma(e)=e,\\
    -1, & \text{if }\sigma(e)\neq e.
  \end{cases}
\end{equation}
The resulting sign assignment is then grouped
by link component to match the inversion datum format used throughout the
package.

Candidate sign assignments are tested in parallel using Python
\texttt{multiprocessing}.
For each candidate, the local validation step computes the crossing sign matrix
and rejects assignments that do not realize one of the allowed local patterns
\eqref{eq:positive-inversion-patterns}--\eqref{eq:negative-inversion-patterns} as a sanity check.
Candidates passing this filter are then tested by the Phase~2 constraint machinery.

\subsection{Phase 2: Linear Constraint System and ILP Export}\label{sec:phase-2}

Starting from a braid and an inversion datum, the code generates the
following families of relations.

\paragraph{Boundary conditions}
The reduced trace fixes the endpoints of the segment~0 to the ground
state. In the code, this is the pair of constraints
\begin{equation}\label{eq:boundary}
  \mathrm{state}(0,0)\equiv \mathrm{state}(0,n)=s_0(b_1),
\end{equation}
with \(s_0(b_1)=0\) when the incoming segment of the open strand $b_1$ is marked \(+\) and \(s_0(b_1)=-1\) when it is
marked \(-\) (cf.~\eqref{eq:ground-state}).

\paragraph{Closure (periodicity)}
For every strand index \(i\), the closure imposes the the equivalences
\eqref{eq:periodicity}.

\paragraph{Sign bounds}
The inversion datum specifies whether $\mathrm{state}(i,t)$ takes values in the nonnegative or
negative integers, for $i=0,\dots, N-1$ and $t=0,\dots, n$.
In the code, a sign assignment is represented as a function on the lattice locations \((i,t)\mapsto \pm 1\).
For each state variable \(u\) (associated to a given segment), this yields one inequality of the
form
\begin{equation}\label{eq:sign-bounds}
  u\ge 0\quad\text{or}\quad u\le -1.
\end{equation}
These are stored as \texttt{Leq} constraints against the literals
\texttt{[0]} and \texttt{[-1]}.

\paragraph{Crossing relations}
At every crossing, the conservation law \eqref{eq:state-conservation} is added as a relation.
In addition, for crossings of type \(R1\) or \(R4\), the code adds two inequalities
ensuring that the $q$-binomial and $q$-Pochhammer indices are nonnegative.
In the coordinate convention \eqref{eq:corner-coordinates}, these are
\begin{equation}\label{eq:r1-r4-ineq}
  \begin{aligned}
  R1:&\qquad \mathrm{state}(k,t+1)\le \mathrm{state}(k-1,t),\\
     &\qquad \mathrm{state}(k,t)\le \mathrm{state}(k-1,t+1),\\
  R4:&\qquad \mathrm{state}(k-1,t)\le \mathrm{state}(k,t+1),\\
     &\qquad \mathrm{state}(k-1,t+1)\le \mathrm{state}(k,t).
  \end{aligned}
\end{equation}
For types \(R2\) and \(R3\) the corresponding nonnegativity conditions are
automatic from the sign bounds \eqref{eq:sign-bounds} (the relevant differences
are strictly positive because the two diagonal entries have opposite signs), so
no extra inequalities are required.

\paragraph{Reduction and parametrization}
The raw relations are simplified by a fixed-point reduction:
constants are propagated through braid closures, symmetric inequalities \((a\le b\) and
\(b\le a)\) are turned into equalities, and vacuous constraints are removed.
From the reduced system, the code assigns fresh
symbols to the remaining free variables and propagates these to all state
variables using the equalities and closure equivalences.
The outcome is an explicit affine-linear expression
\begin{equation}\label{eq:state-affine}
  \mathrm{state}(i,t)=c_{i,t}+\sum_{r=1}^m M_{i,t,r}\,a_r,
\end{equation}
in a small set of free integer parameters \(a_1,\dots,a_m\).

Finally, the code rewrites all inequalities and all $x$-degree truncation
conditions in these parameters.
Single-variable inequalities are used to normalize the sign of each parameter:
if a parameter is constrained by \(a_r\le -1\), the substitution
\(a_r=-b_r-1\) is performed so that all parameters become nonnegative.

At this stage, the global boundedness check is performed. After adding the degree truncation inequalities, the HiGHS-based routine checks that the resulting integer polyhedron, defined by the reduced system, is nonempty and bounded. This is done by maximizing each free integer parameter in turn and rejecting the candidate if any direction is unbounded.
If successful, the resulting data---the braid, the components of its closure, a list of linear
inequalities \(0\le \alpha+\sum \beta_r b_r\), and the full assignment
matrix \eqref{eq:state-affine} evaluated at every lattice location---is
serialized into a compact CSV file.
This file is the only input consumed by the C\texttt{++} backend.

\subsection{Phase 3: Integer-Point Enumeration and \texorpdfstring{$R$-Matrix}{R-Matrix} Evaluation}\label{sec:phase-3}

Phase~3 is performed by the executable \texttt{fk\_main}.
It parses the CSV into an \texttt{FKConfiguration} consisting of:
the braid word (generator indices), the crossing types \(1,2,3,4\) corresponding
to \(R1,R2,R3,R4\), the component indices for each crossing and each closed
strand, a list of linear ``criteria'' rows encoding the $x$-degree truncation,
auxiliary inequalities from Phase~2, and the affine assignment matrix
\eqref{eq:state-affine} at all lattice locations.

\paragraph{Validated bounded criteria}
The C\texttt{++} code must turn the degree constraints into effective upper
bounds for all free parameters.
Because the criteria are stored as general linear forms in the nonnegative
parameters, they do not necessarily exhibit an explicit upper bound for each
parameter.
The routine \texttt{findValidCriteria} therefore searches for an \emph{auxiliary}
collection of bounding rows, in which every parameter appears with a negative
coefficient in some row whose remaining coefficients are nonpositive.
If the original degree criteria already have this monotonicity property, they
are used directly.
Otherwise, \texttt{fk\_main} performs a breadth-first exploration in the space of
rows obtained by repeatedly replacing a criterion row \(c\) by \(c+\lambda u\), with \(\lambda=\tfrac12\) in the current implementation, where \(u\) ranges over the supporting inequalities\footnote{The choice \(\lambda=\tfrac12\) is an implementation
  heuristic rather than a
  mathematical requirement; any nonnegative \(\lambda\)
  would preserve redundancy.}.
Since every supporting inequality is itself enforced during enumeration,
replacements of the form \(c\mapsto c+\lambda u\) only weaken the row constraint
and are therefore redundant; they are used solely to derive explicit upper bounds
for the parameters.

\paragraph{Enumerating integer points}
After a bounded criterion set is found, the admissible parameter vectors are
precisely the integer points
\begin{equation}\label{eq:polytope}
  \mathcal P=
  \Bigl\{b\in\mathbb Z_{\ge 0}^m:\;\text{all exported inequalities hold}\Bigr\}.
\end{equation}
The enumerator first derives per-variable bounds from the validated criteria and then performs an iterative depth-first search over variable assignments.
At each step, a chosen inequality with a negative coefficient in the next variable provides an explicit upper bound for that variable after substituting the already-fixed values.
The first variable is parallelized with OpenMP, and each branch checks the full inequality set before accepting a point.
The resulting list of points is exactly the set of admissible integer parameter assignments for which the truncated inverted sum is finite.

\paragraph{Evaluating the $R$-matrix product}
For each point \(b\in\mathcal P\), the engine evaluates the affine assignment
matrix to recover the numerical state values at every lattice location
\((i,t)\); in particular, at each crossing one obtains the quadruple
\((i,j,i',j')\) used in
Section~\ref{sec:large-color-r-matrices}.
The contribution of \(b\) to the inverted state sum is computed in a factorized form by \eqref{eq:positive-rmatrix-1}--\eqref{eq:negative-rmatrix-1}:
the monomial prefactors in \(q\) and the \(x_c\) are accumulated as explicit
exponent offsets, while the $q$-binomial and (inverse) $q$-Pochhammer factors are
expanded as polynomials and multiplied in.
Concretely, for each crossing type \(r\in\{1,2,3,4\}\), the code multiplies a $q$-binomial
\(\begin{bmatrix}\cdot\\\cdot\end{bmatrix}_q\) with either \((xq^{\ast};q)_{\ast}\) or its inverse.
Inverse Pochhammer factors are expanded as geometric series at $x_i=0$ for all $i=1,\dots,\ell$ and truncated using
the remaining $x$-degree budget so that no term beyond the requested degree is computed.
Polynomial arithmetic is performed using FLINT multivariate polynomials
together with thread-safe caches for $q$-binomials,
Pochhammers, and per-crossing products.

\paragraph{Normalization and output}
After summing the contributions over all points, the final reduced-trace
normalization is applied by multiplying by \((-1+x_1)\) and truncating in each topological variable.
The result is written to JSON together with metadata recording the number of
variables and the constant fractional offsets in the \(q\)- and \(x\)-powers.
The Python layer reads this JSON and, if requested, converts it to a symbolic
expression by re-introducing the fractional power offsets.

\section{Usage and examples}\label{sec:usage-examples}

The package is distributed as a Python library, a command-line tool, and a
Mathematica paclet. All three interfaces ultimately call the same pipeline
from Section~\ref{sec:algorithm}:
Phase~1 (multicycle enumeration), Phase~2 (boundedness check and ILP export), and
Phase~3 (compiled enumeration and $R$-matrix evaluation). In this section we
describe the supported entry points and illustrate typical workflows.

\subsection{Braid input conventions}\label{sec:braid-input-conventions}

Throughout the package a braid on $N$ strands is encoded as a list of signed
generator indices
\[
  \beta=[g_1,\dots,g_n],\qquad g_t\in \pm\{1,\dots,N-1\}.
\]
The integer $|g_t|=k$ corresponds to the Artin generator $\sigma_k$ and the sign
of $g_t$ specifies the crossing sign (positive for $\sigma_k$ and negative for
$\sigma_k^{-1}$). The computation gives the invariant associated to the {closure} of this braid word.

On the command line, braid words may be passed as a string in any of the
following equivalent formats:
\begin{verbatim}
"[1, -2, 1, -2]"     (JSON-style)
"1,-2,1,-2"          (comma-separated)
"1 -2 1 -2"          (space-separated)
\end{verbatim}
In a shell, the braid argument should be quoted so that a negative generator
is not misinterpreted as a command-line flag.

\subsection{Command-line interface}\label{sec:command-line-interface}

After installation, the entry point \texttt{fk} is available.
Running \texttt{fk --help} shows the global help as well as the available
subcommands. The CLI is intentionally split into a minimal ``one-off'' mode
for quick exploration and a configuration-file mode which exposes the full set
of computation parameters.

\paragraph{Interactive wizard}
Invoking \texttt{fk} with no additional arguments launches an interactive
terminal wizard (a richer version is used when the optional dependency
\texttt{rich} is installed). The wizard prompts for a braid word and degree,
offers a preset choice (single-threaded vs. parallel), and then runs the
three-phase pipeline while reporting progress.

In addition to the default wizard, one may explicitly run
\begin{verbatim}
fk interactive
fk interactive --enhanced
fk interactive --quick
\end{verbatim}
The \texttt{--quick} option skips nonessential prompts and uses sensible
defaults. In either mode the wizard can optionally (i) save intermediate files
to disk (inversion JSON, ILP CSV, and the final result JSON), and (ii) print a
symbolic rendering of the output using SymPy (if installed).

\paragraph{Simple mode (single braid)}
The subcommand
\begin{verbatim}
fk simple "[1,1,1]" 2
\end{verbatim}
computes the truncated invariant for the closure of the braid and prints the
result as JSON. This ``simple'' mode is a convenience wrapper intended for
quick experiments; it uses defaults for parallelism and saving, and exposes
only a small number of options.

If SymPy is installed (\texttt{pip install "fkcompute[symbolic]"}), one may ask
for a human-readable expression instead of raw JSON:
\begin{verbatim}
fk simple "[1,1,1]" 2 --symbolic
fk simple "[1,1,1]" 2 --format latex
fk simple "[1,1,1]" 2 --format mathematica
\end{verbatim}
The \texttt{--format} flag automatically enables symbolic printing and selects
one of \texttt{pretty}, \texttt{inline}, \texttt{latex}, or \texttt{mathematica}.

\paragraph{Reformatting saved output}
When results are saved to disk (see below), the subcommand
\begin{verbatim}
fk print-as result.json --format latex
\end{verbatim}
reformats a previously computed JSON result into a symbolic expression without
rerunning the three-phase pipeline.

\paragraph{Config mode (reproducible runs and batches)}
For reproducibility and large computations we recommend using configuration
files. A template may be created by
\begin{verbatim}
fk template create my_run.yaml
fk template create my_run.yaml --overwrite
\end{verbatim}
The resulting YAML file contains commented documentation for the supported
keys. After editing, run
\begin{verbatim}
fk config my_run.yaml
\end{verbatim}
Configuration files can be either JSON (always available) or YAML (requires
\texttt{PyYAML}, installed via \texttt{pip install "fkcompute[yaml]"}).

The config loader supports two shapes:
\begin{itemize}
  \item \emph{single computation:} a file with top-level keys \texttt{braid} and
    \texttt{degree};
  \item \emph{batch mode:} a file with a top-level \texttt{computations} list,
    each entry containing a \texttt{name}, \texttt{braid}, and \texttt{degree}.
\end{itemize}
In batch mode, any additional top-level keys are treated as defaults and are
merged into each computation unless overridden locally.

\subsection{Output format}\label{sec:output-format}

All entry points return the same JSON-serializable dictionary.
For a single computation, the top-level structure is
\begin{verbatim}
{
  "terms": [...],
  "metadata": {...}
}
\end{verbatim}
The list \texttt{terms} encodes a sparse polynomial in the topological
variables and $q$. Each element has an $x$-exponent vector and a sparse
$q$-polynomial coefficient:
\begin{verbatim}
{
  "x": [a1, a2, ...],
  "q_terms": [
    {"q": k1, "c": "n1"},
    {"q": k2, "c": "n2"},
    ...
  ]
}
\end{verbatim}
Here \texttt{c} is stored as a decimal string to preserve arbitrary-precision
integers produced by the FLINT backend.
The \texttt{metadata} field records auxiliary information such as the braid word
actually used (cyclic shifts may be applied during inversion), the inversion
datum, the number of link components, and fixed fractional power offsets in the
$x$- and $q$-exponents.  When symbolic output is requested, these offsets are
reintroduced automatically.

To reconstruct the series directly from the JSON, let \(C\) be
\texttt{metadata["num\_x\_variables"]}.  For a term with
\texttt{"x": [a1,\dots,aC]}, define
\[
  Q_{\mathbf a}(q)
  =
  \sum_{\{\texttt{"q"}=k,\ \texttt{"c"}=n\}\in\texttt{q\_terms}}
  n q^k,
\]
where the coefficient string \(n\) is first converted to an integer.  Ignoring
the optional overall offsets for the moment, the JSON term contributes
\[
  x_1^{a_1}\cdots x_\ell^{a_\ell} Q_{\mathbf a}(q).
\]
Thus the full truncated output is
\begin{multline}\label{eq:json-to-series}
  F(\mathbf x,q)
  =
  q^{q_0}x_1^{r_1}\cdots x_\ell^{r_\ell} \\
  \times \sum_{\texttt{term}\in\texttt{terms}}
  x_1^{a_1}\cdots x_\ell^{a_\ell}
  \\
  \times\sum_{\texttt{qt}\in\texttt{term["q\_terms"]}}
  c_{\texttt{qt}} q^{k_{\texttt{qt}}},
\end{multline}
where \(q_0=\texttt{metadata["overall\_q\_power"]}\) and
\((r_1,\dots,r_\ell)=\texttt{metadata["overall\_x\_powers"]}\), with missing
offset fields interpreted as zero.  These offsets may be fractional; the
integer exponents stored in \texttt{terms} are the sparse polynomial part after
factoring out the common monomial \(q^{q_0}\mathbf x^{\mathbf r}\).

For example, the JSON fragment
\begin{verbatim}
{
  "terms": [
    {
      "x": [0],
      "q_terms": [
        {"q": 4, "c": "1"}
      ]
    },
    {
      "x": [2],
      "q_terms": [
        {"q": -2, "c": "1"},
        {"q": 0,  "c": "-1"}
      ]
    }
  ],
  "metadata": {
    "num_x_variables": 1,
    "overall_x_powers": [-0.5],
    "overall_q_power": 0.75
  }
}
\end{verbatim}
has one topological variable \(x\).  The first entry has \(x\)-degree \(0\)
and \(q\)-coefficient \(q^4\), the second has \(x\)-degree \(2\) and \(q\)-coefficient \(q^{-2}-1\), and there is an overall shift of \(x\)-degree \(-1/2\)
and \(q\)-degree \(3/4\).  Therefore it represents
\[
  F(x,q)=q^{3/4}x^{-1/2}\left(q^4+x^2(q^{-2}-1)\right).
\]
For a link, the vector \texttt{"x": [a1,a2,\dots]} is read in the same way as
\(x_1^{a_1}x_2^{a_2}\cdots\); for two components the symbolic formatter names
these variables \(x\) and \(y\).

In batch mode (a config file with a top-level \texttt{computations} list), the
return value is a dictionary mapping computation names to per-computation result
dictionaries of the above form.

\subsection{Python API}\label{sec:python-api}

Programmatic access is provided by the single function \texttt{fkcompute.fk}.
In ``simple'' mode, it is called with a braid list and a degree:
\begin{verbatim}
from fkcompute import fk

result = fk([1, 1, 1], degree=2)
\end{verbatim}
The returned dictionary has the same \texttt{terms}/\texttt{metadata} structure as
the CLI output.

The Python entry point additionally accepts keyword arguments controlling the
pipeline. The most commonly used are:
\begin{itemize}
  \item \texttt{threads}: number of threads used by the Phase~3 C\texttt{++} backend;
  \item \texttt{max\_workers} and \texttt{chunk\_size}: parallelism settings for
    the Phase~2 inversion search;
  \item \texttt{save\_data}, \texttt{save\_dir}, \texttt{name}: enable saving of
    the inversion JSON, ILP CSV, and final result JSON;
  \item \texttt{symbolic}: attach a SymPy pretty-printed expression as
    \texttt{metadata["symbolic"]} (requires SymPy).
\end{itemize}

Two additional mechanisms are provided for resuming or accelerating
computations:
\begin{itemize}
  \item passing \texttt{inversion} or \texttt{inversion\_file} skips Phase~1 by
    supplying a sign assignment directly;
  \item passing \texttt{ilp\_file} skips Phases~1--2 and runs Phase~3 on a
    pre-generated ILP export.
\end{itemize}
These options are particularly useful when computing the same braid at several
degrees: Phase~1 can be cached once and reused.

Finally, \texttt{fk} also accepts a configuration path as its first argument.
In that case it dispatches to the same config loader used by \texttt{fk config}.

\subsection{Mathematica interface}\label{sec:mathematica-interface}

The directory \texttt{mathematica/FkCompute/} contains a Wolfram Language paclet
which calls the Python implementation through a JSON bridge
\texttt{python -m fkcompute.mathematica\_bridge}.  After installing the paclet,
one loads it in Mathematica via
\begin{verbatim}
PacletInstall[".../mathematica/FkCompute"]
Needs["FkCompute`"]
\end{verbatim}
and computes (for example) the trefoil by
\begin{verbatim}
FkCompute[{1, 1, 1}, 2]
\end{verbatim}
When SymPy is available in the Python environment, the default return value is
a Mathematica expression for the truncated polynomial in the variables
corresponding to the link components and $q$.

The wrapper also exposes \texttt{FkComputeVersion[]} to query the Python and
\texttt{fkcompute} versions seen by Mathematica.  If Mathematica does not find
the correct Python executable (for instance when using a virtual environment),
set the executable explicitly, e.g.
\begin{verbatim}
SetFkComputePythonExecutable[
  "/path/to/python3"
]
\end{verbatim}

The function \texttt{FkCompute} accepts options mirroring the Python keyword
arguments (e.g. \texttt{"Threads"}, \texttt{"MaxWorkers"}, \texttt{"Verbose"},
\texttt{"SaveData"}, and \texttt{"Preset"}), which are passed through to the
bridge unchanged.

\subsection{Worked examples}\label{sec:worked-examples}

We conclude with three representative examples.

\paragraph{Figure eight knot}
The figure eight knot is the closure of $\sigma_1\sigma_2^{-1}\sigma_1\sigma_2^{-1}$, encoded by \texttt{[1,-2,1,-2]}.
In simple CLI mode:
\begin{verbatim}
fk simple "[1,-2,1,-2]" 10 --symbolic
\end{verbatim}
This example is a homogeneous braid, so Phase~1 is instantaneous.

\paragraph{Knot $8_{20}$}
The braid
$\sigma_{1}^{3}\sigma_{2}^{-1}\sigma_1^{-3}\sigma_{2}^{-1}$ of the knot $8_{20}$ is not homogeneous.
A typical workflow is to use the wizard (to tune workers/threads) or to store a
reproducible configuration:
\begin{verbatim}
# config_820.yaml
braid: [1,1,1,-2,-1,-1,-1,-2]
degree: 5
max_workers: 8
threads: 8
save_data: true
name: k8_20
\end{verbatim}
and run \texttt{fk config config\_820.yaml}.  The saved
\texttt{k8\_20\_inversion.json} can then be reused at higher degrees.

\paragraph{A two-component link (the $(2,4)$ torus link)}
The $(2,4)$ torus link is the closure of $\sigma_1^4$, encoded by \texttt{[1,1,1,1]}.
The computation returns a polynomial in two topological variables and $q$.
For example:
\begin{verbatim}
fk simple "[1,1,1,1]" 2 --format latex
\end{verbatim}
prints a LaTeX expression in variables \texttt{x} and \texttt{y} (one per
component), while the JSON output encodes the same data via $x$-degree vectors
of length two.

\section{Applications}\label{sec:applications}

The package \texttt{fkcompute} was used to generate the largest dataset of known $F_K$ and $F_L$ invariants yet, which is available at \cite{topologyfyi}.
The package \texttt{fkcompute} has been used to obtain several observations, some of which are now theorems. Here we list two major examples:

\paragraph{Leading term}

Using $\texttt{fkcompute}$, we observed that the coefficient of the minimal $x$-power is a monomial, which we subsequently proved in \cite{OSSS25}. For all knots up to 12 crossings, the corresponding power of $q$ equals $g(K)-\lambda(K)$, where $g(K)$ is the Seifert genus of the knot $K$, and $\lambda(K)$ is its \emph{Hopf invariant} (introduced in \cite{Rudolph1987} as enhanced Milnor number). This holds for any homogeneous braid knot as well, whereas it remains an interesting conjecture for an arbitrary fibered knot.

\paragraph{Slopes}

Consider the $F_K$ invariant as a series of Laurent polynomials $f_n(q)$, $n = 0, 1,2 \dots$, given by the coefficients of $x^n$. The minimal powers of $f_n(q)$ form a quadratic quasipolynomial, whose leading coefficient, the $F_K$ \emph{slope}, is conjecturally a boundary slope of the knot $K$. Using $\texttt{fkcompute}$, we identified the slopes for fibered knots up to 10 crossings (cf. \cite{OSSS25}).

In particular, we demonstrated that there are knots like $9_{44}$ with slope strictly smaller than $-1/2$, disproving the conjectural meromorphicity of Park's inverted Habiro series \cite{Park21,Svo2025}.

\section{Performance}\label{sec:performance}

The cost of a computation splits naturally along the three phases of
Section~\ref{sec:algorithm}.  In practice, Phase~1 (multicycle enumeration) is
a fixed, degree-independent overhead for a given braid, Phase~2 contains the
degree-dependent boundedness check and ILP export, and Phase~3 (integer point
enumeration and $R$-matrix evaluation) dominates at high truncation degree.
All timings reported below were obtained on a desktop workstation with
a 12-core AMD Ryzen~9 7900X CPU, running the compiled backend with 12 threads
unless stated otherwise.

\paragraph{Phase 1: multicycle enumeration}
The naive search space of sign assignments for a braid word with $n$ crossings
has size $2^{2n}$, since the braid diagram has $2n$ segments.  The
permutation-based enumeration of multicycle candidates described in
Section~\ref{sec:phase-1} reduces this dramatically: across the $1245$
\emph{nice} knots with at most 12 crossings in our test suite, whose braid
presentations contain up to 18 crossings, the number of locally admissible
candidates never exceeds $310$ (median $134$), compared to a generic count of
up to $2^{36}\approx 6.9\times 10^{10}$; the median reduction factor is about
$6.5\times 10^{-6}$ even at moderate braid lengths.  Enumerating the candidates
takes milliseconds; the subsequent HiGHS-based boundedness certification is a
Phase~2 step and is parallelized over candidates.  For the braids in our test
range, the combined candidate-search-and-certification stage typically
completes in well under a second (e.g. 0.13\,s for the non-homogeneous
$8_{20}$ braid of Section~\ref{sec:worked-examples}).  For homogeneous braids,
Phase~1 is a deterministic rule and its cost is negligible.

\paragraph{Phase 3: scaling with the truncation degree}
The number of integer points in the polyhedron \eqref{eq:polytope} grows
rapidly with the requested $x$-degree, and this growth dictates the overall
scaling.  Table~\ref{tab:timings} collects representative wall-clock times.
Low-complexity knots are essentially instantaneous even at high degree: the
figure-eight knot reaches degree 50 in under two seconds.  For braids with more
crossings the degree becomes the limiting factor; for the 12-crossing braid of
$12n_{242}$, each increment of the degree by 5 increases the run time roughly
by an order of magnitude beyond degree 20.
Links are comparable to knots of
similar diagrammatic complexity, with the truncation applied per component
variable.

\begin{table}
\centering
\scriptsize
\setlength{\tabcolsep}{4pt}
\begin{tabular}{l r r r}
\hline
input & braid & degree & time (s) \\
\hline
$3_1$ & \texttt{[1,1,1]} & 10 & $<0.1$ \\
$4_1$ & \texttt{[1,-2,1,-2]} & 30 & $0.1$ \\
$4_1$ & \texttt{[1,-2,1,-2]} & 50 & $1.6$ \\
$8_{20}$ & \texttt{[1,1,1,-2,-1,-1,-1,-2]} & 5 & $0.2$ \\
$11a_{128}$ & \texttt{[1,-2,-1,3,-2,-1,3,-2,-4,3,5,-4,5]} & 10 & $0.6$ \\
L5a1 & \texttt{[-1,2,-1,2,-1]} & 8 & $0.1$ \\
L8n6 & \texttt{[2,-3,4,-3,-4,2,-3,-4,1,-2,1]} & 8 & $0.5$ \\
$12n_{242}$ & \texttt{[1,2,2,1,1,2,2,2,2,2,2,2]} & 15 & $0.1$ \\
$12n_{242}$ & \texttt{[1,2,2,1,1,2,2,2,2,2,2,2]} & 20 & $1.1$ \\
$12n_{242}$ & \texttt{[1,2,2,1,1,2,2,2,2,2,2,2]} & 25 & $13$ \\
$12n_{242}$ & \texttt{[1,2,2,1,1,2,2,2,2,2,2,2]}& 30 & $116$ \\
\hline
\end{tabular}
\caption{Representative wall-clock times (12 threads, AMD Ryzen~9 7900X).
Unless noted, an inversion datum was supplied or found instantaneously, so the
times are dominated by Phase~3.}
\label{tab:timings}
\end{table}

\paragraph{Parallel scaling}
The enumeration and evaluation loop is parallelized with OpenMP over the first
free parameter.  On the $12n_{242}$ example at degree 25, the computation takes
$129$\,s on a single thread, $35$\,s on 4 threads, and $13$\,s on 12 threads,
i.e.\ the speedup is close to linear in the number of cores.  Memorization of
$q$-binomials, $q$-Pochhammer factors, and per-crossing partial products in
thread-safe caches contributes substantially to throughput, since the same
local factors recur across many states.

\section{Validation}\label{sec:validation}

The correctness of the implementation is monitored by an automated test suite
distributed with the source code.  It combines external cross-checks against
previously computed data with the theoretical consistency conditions of
Section~\ref{sec:fk-fl-invariants}.

\paragraph{Relations with other quantum invariants}

We use the relation between $F_K$ and the Melvin--Morton--Rozansky (MMR) expansion to provide basic consistency checks for the computations implemented in \texttt{fkcompute}.

For a knot $K$, define
\begin{equation}
    f_K(x,q)=\frac{F_K(x,q)}{x^{\frac12}-x^{-\frac12}}.
\end{equation}
The expansion of $f_K(x,e^\hbar)$ at $\hbar=0$ recovers the Melvin--Morton--Rozansky expansion \cite{BNG,RozMMR} in \eqref{eq:FK-mmr-background}, with the rational functions in $x$ expanded at $x=0$ \cite[Theorem 2]{Park21}.
Consequently,
\begin{equation}\label{eq:q-one-limit-background}
  f_K(x,1)
  =
  \frac{1}{\Delta_K(x)}
\end{equation}
and
\begin{equation}\label{eq:q-derivative-background}
    D_q f_K(x,q)\Big|_{q=1}
  =
  \frac{P_{K,1}(x)}{\Delta_K(x)^3}~,
\end{equation}
where $D_q$ denotes the first-derivative in the $q$-direction.

The specializations \eqref{eq:q-one-limit-background} and
\eqref{eq:q-derivative-background} provide two nontrivial analytic checks.  The
test suite contains an independent pure-Python implementation of one of the
Bar-Natan--van der Veen polynomial-time invariants \cite{BNVdV22}, computed directly from a planar diagram and without any reference
to the $R$-matrix construction. These give the Alexander polynomial $\Delta_K(x)$ and a perturbed Alexander polynomial $\rho_{1,K}(x)$, which conjecturally agrees with the one-loop
polynomial $P_{K,1}(x)$ in \eqref{eq:q-derivative-background}. For each computed example, the
specialization $q=1$ of $f_K(x,q)$ is checked against the
expansion of $1/\Delta_K(x)$, and the first $q$-derivative at $q=1$ is checked
against the expansion of $\rho_{1,K}(x)/\Delta_K(x)^3$, in both cases up to the
truncation degree.

For a link $L$ with more than one component, we compare the evaluation of $F_L(x_1,\dots,x_\ell,q)$ at $q=1$ with the iterated Laurent series of the inverse of its multivariable Alexander polynomial; see \cite[Chapter 2]{Xin04} for a definition. In particular,
\begin{equation}
    \Delta_L(x_1,\dots,x_\ell)\cdot F_L(x_1,\dots,x_\ell,1) = 1.
\end{equation}
Furthermore, we compare the evaluation of the first $q$-derivative at $q=1$ with a multivariable generalization of the perturbed Alexander polynomial $\rho_{1,L}$, defined for a subclass of links including nice links. More details to appear in \cite{PassaroSuarez}.

\paragraph{Internal consistency}
Several invariance properties are tested directly.  The series produced from
different acceptable inversion data for the same braid, and from braid variants
related by cyclic rotation or by the horizontal mirror used in the Phase~2
sign diagram validation, agree on all examples tested, consistent with the fact that the
inverted state sum is independent of these choices \cite{Park21}.
The combinatorial layer is further pinned down by snapshot tests of the Phase~2
constraint export and by tests verifying the permutation-to-sign-diagram
enumeration of Section~\ref{sec:phase-1} against stored reference data for
knots and links.

\section{Scope and limitations}\label{sec:limitations-future-work}

\paragraph{Scope of the method}
The package computes inverted state sums, and therefore applies only when an
acceptable inversion datum exists for some braid variant of the input.  This
covers homogeneous braids and, conjecturally, all fibered knots; but the
existence of such a datum is not guaranteed in general, and the geometric
meaning of the inversion datum remains poorly understood.  When the search
fails at a given degree, the package reports failure rather than falling back
to a different construction; finding either a characterization of \emph{nice}
presentations or a complementary algorithm for the remaining cases is an open
problem with direct practical impact.  Relatedly, the result of Phase~1 depends
on the braid presentation: a braid whose search fails may admit a conjugate or
Markov-equivalent presentation for which it succeeds, and the package currently
explores only cyclic rotations and an optional mirror rather than the full
equivalence class.

\paragraph{Computational scaling}
The number of admissible states grows quickly with both the truncation degree
and the diagrammatic complexity of the braid, and the enumeration cost grows
accordingly (cf.\ Table~\ref{tab:timings}).  This places high-degree
computations for long braids out of reach.  Natural improvements include
smarter pruning of the integer point enumeration, exploiting symmetries of the
constraint polyhedron, and incremental computation of consecutive
$x$-coefficients, which would allow a single run to extend an existing
truncation rather than recompute it.

\section*{Acknowledgements}
The authors are grateful to Sergei Gukov for his support completing this project.
The work of Davide Passaro is supported by a Sherman Fairchild Postdoctoral Fellowship sponsored
by the Walter Burke Institute for Theoretical Physics.
Josef Svoboda was supported by the Simons Collaboration grant New Structures in Low-Dimensional Topology. Lara San Martín Suárez received the support of a fellowship from ``la Caixa'' Foundation (ID 100010434), with fellowship code LCF/BQ/EU23/12010094.

\bibliography{refs}

@article{GM,
  author = "Gukov, Sergei and Manolescu, Ciprian",
  title = "{A two-variable series for knot complements}",
  doi = "10.4171/qt/145",
  journal = {Quantum Topol.},
  fjournal = {Quantum Topology},
  volume = "12",
  number = "1",
  pages = "1--109",
  year = "2021"
}

@misc{OSSS25,
  doi = {10.48550/ARXIV.2512.23700},
  url = {https://arxiv.org/abs/2512.23700},
  author = {Orland,  Paul and San Martín Suárez, Lara and Saunders-A'Court,  Toby and Svoboda,  Josef},
  title = {Quantum Invariants and Fiberedness},
  publisher = {arXiv},
  year = {2025},
  copyright = {arXiv.org perpetual,  non-exclusive license}
}

@article{huangfu2018parallelizing,
  title={Parallelizing the dual revised simplex method},
  author={Huangfu, QHJAJ and Hall, JA Julian},
  journal={Mathematical Programming Computation},
  volume={10},
  number={1},
  pages={119--142},
  year={2018},
  publisher={Springer}
}

@article{Rudolph1987,
  title = {{Isolated critical points of mappings from R4 to R2 and a natural splitting of the {Milnor} number of a classical fibered link. Part I: Basic theory; examples}},
  volume = {62},
  ISSN = {1420-8946},
  url = {http://dx.doi.org/10.1007/BF02564467},
  DOI = {10.1007/bf02564467},
  number = {1},
  journal = {Commentarii Mathematici Helvetici},
  publisher = {European Mathematical Society - EMS - Publishing House GmbH},
  author = {Rudolph,  Lee},
  year = {1987},
  month = dec,
  pages = {630--645}
}

@misc{Svo2025,
  title={Inverted {Habiro} Series and its Residues}, 
  author={Josef Svoboda},
  year={2025},
  eprint={2509.22610},
  archivePrefix={arXiv},
  primaryClass={math.GT},
  url={https://arxiv.org/abs/2509.22610}, 
}

@phdthesis{ParkThesis,
  doi = {10.7907/M1FR-6038},
  author = {Park,  Sunghyuk},
  school = {California Institute of Technology},
  title = {3-Manifolds,  Q-Series,  and Topological Strings},
  publisher = {California Institute of Technology},
  year = {2022},
  copyright = {No commercial reproduction,  distribution,  display or performance rights in this work are provided.}
}

@article{MM,
  doi = {10.1007/bf02099310},
  url = {https://doi.org/10.1007/bf02099310},
  year = {1995},
  month = may,
  publisher = {Springer Science and Business Media {LLC}},
  volume = {169},
  number = {3},
  pages = {501--520},
  author = {P. M. Melvin and H. R. Morton},
  title = "{The coloured Jones function}",
  journal = {Comm. Math. Phys.}
}

@article{RozMMRConj,
	doi = {10.1007/bf02506408},
	url = {https://doi.org/10.1007\%2Fbf02506408},
	year = 1997,
	month = jan,
	publisher = {Springer Science and Business Media {LLC}
},
	volume = {183},
	number = {2},
	pages = {291--306},
	author = {L. Rozansky},
	title = "{Higher order terms in the {Melvin--Morton} expansion of the colored {Jones} polynomial}",
	journal = {Comm. Math. Phys.}
}

@article{BNG,
  doi = {10.1007/s002220050070},
  url = {https://doi.org/10.1007/s002220050070},
  year = {1996},
  month = may,
  publisher = {Springer Science and Business Media {LLC}},
  volume = {125},
  number = {1},
  pages = {103--133},
  author = {Dror Bar-Natan and Stavros Garoufalidis},
  title = "{On the Melvin--Morton--Rozansky conjecture}",
  journal = {Inventiones Mathematicae}
}

@article{GPV,
  author  = {Gukov, Sergei and Putrov, Pavel and Vafa, Cumrun},
  title   = {Fivebranes and 3-manifold homology},
  journal = {Journal of High Energy Physics},
  year    = {2017},
  volume  = {2017},
  number  = {7},
  pages   = {71},
  doi     = {10.1007/JHEP07(2017)071},
  url     = {https://doi.org/10.1007/JHEP07(2017)071}
}

@article{GPPV,
  doi = {10.1142/s0218216520400039},
  url = {https://doi.org/10.1142/s0218216520400039},
  year = {2020},
  month = feb,
  publisher = {World Scientific Pub Co Pte Lt},
  volume = {29},
  number = {02},
  pages = {2040003},
  author = {Sergei Gukov and Du Pei and Pavel Putrov and Cumrun Vafa},
  title = "{{BPS} spectra and 3-manifold invariants}",
  journal = {J. Knot Th. Ram.},
  fjournal = {Journal of Knot Theory and Its Ramifications},
}

@article{Kh,
  author = {Mikhail Khovanov},
  title = {A categorification of the {Jones} polynomial},
  journal = {Duke Mathematical Journal},
  year = {2000},
  volume = {101},
  number = {3},
  pages = {359--426},
  doi = {10.1215/S0012-7094-00-10131-7},
  publisher = {Duke University Press}
}

@article{Witten,
  author={Edward Witten},
  title={Quantum field theory and the {Jones} polynomial},
  journal={Comm. Math. Phys.},
  year={1989},
  volume={121},
  number={3},
  pages={351--399},
  publisher={Springer}
}

@article{RT,
  author = {Reshetikhin, N. and Turaev, V. G.},
  title = {Invariants of 3-manifolds via link polynomials and quantum groups},
  journal = {Inventiones mathematicae},
  year = {1991},
  volume = {103},
  number = {1},
  pages = {547--597},
  doi = {10.1007/BF01239527}
}

@article{RozMMR,
  title = {The Universal {R-Matrix},  {Burau} Representation,  and the {Melvin--Morton} Expansion of the Colored {Jones} Polynomial},
  volume = {134},
  ISSN = {0001-8708},
  url = {http://dx.doi.org/10.1006/aima.1997.1661},
  DOI = {10.1006/aima.1997.1661},
  number = {1},
  journal = {Adv. in Math.},
  publisher = {Elsevier BV},
  author = {Rozansky,  L.},
  year = {1998},
  month = mar,
  pages = {1--31}
}

@misc{Park21,
  title={Inverted state sums, inverted {Habiro} series, and indefinite theta functions}, 
  author={Sunghyuk Park},
  year={2021},
  eprint={2106.03942},
  archivePrefix={arXiv},
  primaryClass={math.GT},
  url={https://arxiv.org/abs/2106.03942}, 
}

@misc{PassaroSuarez,
    author={Passaro, Davide and San Martín Suárez, Lara},
    note   = {In preparation},
    }

@article{BNVdV22,
  title = {A {perturbed-Alexander} invariant},
  volume = {15},
  ISSN = {1664-073X},
  url = {http://dx.doi.org/10.4171/QT/206},
  DOI = {10.4171/qt/206},
  number = {3},
  journal = {Quantum Topology},
  publisher = {European Mathematical Society - EMS - Publishing House GmbH},
  author = {Bar-Natan,  Dror and van der Veen,  Roland},
  year = {2024},
  month = Apr,
  pages = {449–472}
}

@misc{fkcompute,
  author       = {Orland, Paul and Passaro, Davide and San Mart\'in Su\'arez, Lara  and Saunders-A'Court, Toby and Svoboda, Josef},
  title        = {{fkcompute}},
  year         = {2026},
  howpublished = {\url{https://github.com/caltech-quantum-topology/fkcompute}},
  note         = {GitHub repository}
}

@misc{topologyfyi,
  title        = {},
  author       = {},
  howpublished = {\url{https://topology.fyi/}},
  year         = {2026}
}

@article{crane1994,
   title={Four-dimensional topological quantum field theory, {Hopf} categories, and the canonical bases},
   volume={35},
   ISSN={1089-7658},
   url={http://dx.doi.org/10.1063/1.530746},
   DOI={10.1063/1.530746},
   number={10},
   journal={Journal of Mathematical Physics},
   publisher={AIP Publishing},
   author={Crane, Louis and Frenkel, Igor B.},
   year={1994},
   month=Oct, pages={5136–5154} }

@phdthesis{Xin04,
  author = {Xin, Guoce},
  title  = {The Ring of {Malcev--Neumann} Series and the Residue Theorem},
  school = {Brandeis University},
  year   = {2004},
  eprint = {math/0405133},
  archivePrefix = {arXiv}
}
\bibliographystyle{ytphys.bst}

\end{document}